\documentclass[12pt]{article}

\usepackage{theorem,amssymb}
 
\advance \topmargin by -\headheight
\advance \topmargin by -\headsep
\evensidemargin \oddsidemargin
\author{J.-P. Allouche \\
CNRS, LRI, B\^atiment 490 \\
F-91405 Orsay Cedex (France) \\
{\tt allouche@lri.fr} \\
\and
M. Mend\`es France \\
Math\'ematiques, Universit\'e Bordeaux I \\
F-33405 Talence Cedex (France) \\
{\tt michel.mendes-france@math.u-bordeaux1.fr} \\
}

\title{Euler, Pisot, Prouhet-Thue-Morse, Wallis and 
the duplication of sines}

\date{ }

\def \proof{\bigbreak\noindent{\it Proof.\ \ }}

\def \endpf{{\ \ $\Box$ \medbreak}}

\newtheorem{theorem}{Theorem}
\newtheorem{lemma}{Lemma}
\newtheorem{corollary}{Corollary}
\newtheorem{proposition}{Proposition}
\theorembodyfont{\rm}

\newtheorem{remark}{Remark}

\begin{document}

\maketitle

\section*{Introduction}

If the product obtained by iterating the duplication of
cosines has a well-known simple closed form, namely
$$
\cos x \cos 2x \cos 4x \cdots \cos 2^n x = \frac{\sin 2^{n+1}x}{2^{n+1} \sin x},
$$
there is no such formula for the analogous product of sines.
Nevertheless the product:
$$
P_n(x) := \sin x \sin 2x \sin 4x \cdots \sin 2^n x 
$$
occurs several times in the literature (see for example
\cite{Gelfond, Newman, Newman-Slater, Eminyan1, Eminyan2, Fouvry-Mauduit}), 
where in particular the quantities
$\Vert P_n \Vert_{\infty} := \sup_{x \in {\mathbb R}} |P_n(x)|$
and 
$\Vert P_n \Vert_1 := \int_0^{\pi} |P_n(x)| dx$
enter the picture.

\bigskip

We will give a general identity, involving the sum of binary digits of
an integer, from which formulas for the product $P_n(x)$ and related 
quantities will be deduced. Exploiting these identities will lead us from 
Pisot numbers to Euler and Wallis formulas, from the Prouhet-Thue-Morse 
sequence to asymptotic formulas in analytic number theory, from multigrade 
equalities to the natural introduction of morphisms of the free monoid on 
two letters. 

\section*{Two propositions}

\subsection*{An easy identity}

The following result is not difficult to prove.

\begin{proposition}
Let $a$ be a complex number. Denote by $s(k)$ 
the sum of the binary digits of the nonnegative integer $k$. 
Then the following identity holds in ${\mathbb C}[[X]]$.
\begin{equation}\label{simple}
\prod_{0 \leq k < n} \left(1+aX^{2^k}\right) 
= \sum_{0 \leq j < 2^n} a^{s(j)} X^j.
\end{equation}
\end{proposition}

\proof This is a direct consequence of the uniqueness of the
base $2$ expansion of the integers in $[0, 2^n)$. Alternatively
this can be proved by induction on $n$. \endpf

\begin{corollary}\label{corsimple}
The following well-known identities hold:
$$
\begin{array}{llll}

\displaystyle\prod_{0 \leq k < n} &\left(1+X^{2^k}\right) 
&=& \displaystyle\frac{1-X^{2^n}}{1-X} \\

\displaystyle\,\,\, \prod_{k \geq 0} &\left(1+X^{2^k}\right) 
&=& \displaystyle\frac{1}{1-X} \\

\displaystyle\prod_{Ã0 \leq k < n} &\left(1-X^{2^k}\right)
&=& \displaystyle\sum_{0 \leq j < 2^n} (-1)^{s(j)} X^j \\

\displaystyle\,\,\, \prod_{k \geq 0} &\left(1-X^{2^k}\right) 
&=& \displaystyle\,\,\, \sum_{j \geq 0} (-1)^{s(j)} X^j \\

\end{array}
$$
\end{corollary}

\proof Take $a=1$ and $a=-1$ in Identity~\ref{simple} gives the first 
and third identities. Letting $n$ tend to infinity in the first and
third identities gives the second and fourth identities. \endpf

\begin{remark}

\ { }

-- The first identity in Corollary~\ref{corsimple} above goes back to
L.~Euler (1707-1783). For an unusual occurrence of the Euler formula, 
the reader might want to look at \cite{Mendes-VdP}.

\medskip

-- The third and fourth identities in Corollary~\ref{corsimple} above 
involve the sequence $((-1)^{s(j)})_{j \in {\mathbb N}}$, which is the 
celebrated Prouhet-Thue-Morse sequence. This sequence first appeared 
in 1851 in a paper by E.~Prouhet (1817-1867) for multigrade equalities 
(see below); it was then studied at the beginning of the 20th century 
by A.~Thue (1863-1922) for a question about repetitions in infinite 
sequences (which lead to the nowadays vivid combinatorics of words); in 
the 20's H.C.M.~Morse (1892-1977) rediscovered it for a construction in 
differential geometry: this and more can be found, e.g., in 
\cite{ubiquitous} and the references therein. Also note that K.~Mahler
(1903-1988) studied this sequence, without calling it the Thue-Morse
sequence: in \cite{Mahler1} he proved that the correlation spectrum of 
this sequence contains a nonzero singular continuous component (Kakutani 
proved in the 60's that the spectrum is singular continuous); Mahler told
the second author that \cite{Mahler1} was the first paper he wrote in 
English -- with the help of N.~Wiener. Speaking of Mahler and the 
Thue-Morse sequence, it is worth mentioning that the transcendence 
results he proved in \cite{Mahler2} directly imply the transcendence 
of the real number whose $b$-ary expansion is the Thue-Morse sequence.

\end{remark}

\subsection*{A generalization}

The following generalization of Identity~\ref{simple} holds; it will
be the core of our paper.

\begin{theorem}
Let $a$ be a complex number. If $j=\sum_q e_q(j) 2^q$
with $e_q(j) \in \{0, 1\}$ is the binary expansion of the nonnegative integer 
$j$, and if $\lambda := (\lambda_j)_{j \in {\mathbb N}}$ is a sequence
of complex numbers, let $u(j, \lambda) := \sum_q e_q(j) \lambda_q$.
Then the following identity holds in ${\mathbb C}[[Z]]$.
\begin{equation}\label{step}
\prod_{0 \leq k < n} \left(1 + a \exp (\lambda_k Z) \right)
= \sum_{0 \leq j < 2^n} a^{s(j)} \exp (u(j, \lambda) Z).
\end{equation}
\end{theorem}

\begin{remark}

\ { }

-- In Identity~\ref{step} the formal series $\exp Y$
is defined as usual by $\exp Y := \sum_{n \geq 0} Y^n/n!$.
This identity could be written in the more visual -- but
maybe less conventional -- form
$$
\prod_{0 \leq k < n} \left(1 + a X^{\lambda_k} \right)
= \sum_{0 \leq j < 2^n} a^{s(j)} X^{u(j, \lambda)}.
$$

\medskip

-- Taking $\lambda_j := 2^j$ and $X := \exp Z$ in 
Identity~\ref{step} clearly implies Identity~\ref{simple}.

\end{remark}

\section*{More product formulas}

The identities above can be turned into trigonometric identities
including in particular a formula for the duplication of sines
and cosines.

\begin{proposition}\label{more}

Let $a$ be a complex number and 
$\lambda := (\lambda_j)_{j \in {\mathbb N}}$ be a sequence
of complex numbers. If the binary expansion of the nonnegative 
integer $j$ is $j = \sum_q e_q(j) 2^q$, with $e_q(j) = 0, 1$, 
define $\sigma_q(j)$ by $\sigma_q(j) := 2 e_q(j) -1$. Then, the
following identity holds true in ${\mathbb C}[[Z]]$.
\begin{equation}\label{general}
a^n \prod_{0 \leq k < n} \left(\exp(\lambda_k Z) 
+ a^{-1} \exp(-\lambda_k Z)\right) =
\sum_{0 \leq j < 2^n} a^{s(j)} 
\exp\left(\left(\sum_{0 \leq q < n} 
\sigma_q(j) \lambda_q\right)Z\right).
\end{equation}
\end{proposition}

\proof  Dividing the two sides of Identity~\ref{step} by
$\exp(\frac{1}{2}(\lambda_0 + \lambda_1 + \cdots + \lambda_{n-1})Z)$
yields
$$
\prod_{0 \leq k < n}
\left(\exp\left(-\frac{1}{2}\lambda_k Z\right) 
+ a\exp\left(\frac{1}{2}\lambda_k Z\right)\right)
= \sum_{0 \leq j < 2^n} a^{s(j)}
\exp\left((u(j, \lambda) - 
\frac{1}{2}(\lambda_0 + \cdots \lambda_{n-1})Z\right).
$$
Since
$$
u(j, \lambda) 
- \frac{1}{2}(\lambda_0 + \cdots \lambda_{n-1})
= \frac{1}{2} \sum_{0 \leq q < n} (2e_q(j)-1)\lambda_q
= \frac{1}{2} \sum_{0 \leq q < n} \sigma_q(j) \lambda_q
$$
we get, after factoring $a$ out and replacing $Z$ by $2Z$,
$$
\prod_{0 \leq k < n}
a^n \left(\exp(\lambda_k Z) + a^{-1} \exp(-\lambda_k Z)\right)
= \sum_{0 \leq j < 2^n} a^{s(j)}
\exp\left(\left(\sum_{0 \leq q < n} \sigma_q(j) \lambda_q\right)Z\right)
$$
which is Identity~\ref{general}.

\bigskip

Identity~\ref{general} implies several product identities.

\begin{corollary}
Let $a$ be a complex number and 
$\lambda := (\lambda_j)_{j \in {\mathbb N}}$ be a sequence
of complex numbers. If the binary expansion of the nonnegative
integer $j$ is $j = \sum_q e_q(j) 2^j$, with $e_q(j) = 0, 1$,
define $\sigma_q(j)$ by $\sigma_q(j) := 2 e_q(j) -1$. Then,
\begin{eqnarray}
\label{cosh}
\prod_{0 \leq k < n} \cosh \lambda_k 
&=& 2^{-n} \sum_{0 \leq j < 2^n} 
\cosh\left(\sum_{0 \leq q < n} \sigma_q(j) \lambda_q\right)
\\
\label{cos}
\prod_{0 \leq k < n} \, \cos \lambda_k 
\,\,\, &=& 2^{-n} \sum_{0 \leq j < 2^n} 
\cos\left(\sum_{0 \leq q < n} \sigma_q(j) \lambda_q\right) \\
\label{sinh}
\prod_{0 \leq k < 2m} \sinh \lambda_k 
&=& 2^{-2m} \sum_{0 \leq j < 2^{2m}} (-1)^{s(j)} 
\cosh\left(\sum_{0 \leq q < 2m} \sigma_q(j) \lambda_q\right) \\
\label{sin}
\prod_{0 \leq k < 2m} \, \sin \lambda_k 
\,\,\, &=& (-1)^m 2^{-2m} \sum_{0 \leq j < 2^{2m}} (-1)^{s(j)} 
\cos\left(\sum_{0 \leq q < 2m} \sigma_q(j) \lambda_q\right)
\end{eqnarray}
Furthermore
\begin{eqnarray}
\label{cos2}
\prod_{0 \leq k < n} \cos(2^k x) &=&
2^{-n} \sum_{0 \leq j < 2^n} \cos((2j+1)x - 2^n x) \\
\label{sin2}
\prod_{0 \leq k < 2m} \sin(2^k x) &=&
(-1)^m2^{-2m} \sum_{0 \leq j < 2^{2m}}
(-1)^{s(j)} \cos((2j+1)x - 2^{2m}x)
\end{eqnarray}
Hence
\begin{eqnarray}
\label{cos2simple}
\prod_{0 \leq k < n} \cos(2^k x) &=&
2^{1-n} \sum_{0 \leq j < 2^{n-1}} \cos((2j+1)x) 
\ \ \left( = \frac{\sin 2^n x}{2^n \sin x} \right) \\
\label{sin2simple}
\prod_{0 \leq k < 2m} \sin(2^k x) &=&
(-1)^{m+1}2^{1-2m} \sum_{0 \leq j < 2^{2m-1}}
(-1)^{s(j)} \cos((2j+1)x)
\end{eqnarray}

\end{corollary}

\proof  Taking $a=1$ in Identity~\ref{general}, we see that 
the left-hand expression is invariant under $Z \to -Z$; hence 
the right-hand expression is also invariant under $Z \to -Z$; 
hence adding Identity~\ref{general} with $Z$ to 
Identity~\ref{general} with $-Z$ yields Identity~\ref{cosh}. 
Changing then $\lambda_k$ into $i\lambda_k$ yields 
Identity~\ref{cos}. 

\bigskip

Taking $a=-1$ in Identity~\ref{general}, we see that changing $Z$ into 
$-Z$ multiplies out the left-hand expression by $(-1)^n$. Taking $n=2m$
gives 
$$
\prod_{0 \leq k < n}
\left(\exp(\lambda_k Z) - \exp(-\lambda_k Z)\right)
= \sum_{0 \leq j < 2^n} (-1)^{s(j)}
\exp\left(\left(\sum_{0 \leq q < n} \sigma_q(j) \lambda_q\right)Z\right)
$$
whose left-hand member (hence right-hand member as well) is invariant under 
$Z \to -Z$. Arguing as above, dividing by $2^{2m}$, and taking $Z=1$  
gives Identity~\ref{sinh}. Replacing $\lambda_k$ by $i\lambda_k$
yields Identity~\ref{sin}.

\bigskip

Then Identities~\ref{cos2} and \ref{sin2} are deduced from 
Identities~\ref{cos} and \ref{sin} respectively by taking 
$\lambda_k := 2^k x$ and by noting that 
$\sum_{0 \leq q < n} \sigma_q(j) 2^j = (2j+1) - 2^n$. 

\bigskip

Finally if we split the sum in Identity~\ref{cos2} as
$\sum_{0 \leq j < 2^{n-1}} + \sum_{2^{n-1} \leq j < 2^n}$ and if we
put $j = 2^{n-1} - t -1$ in the first sum and $j = 2^{n-1} + t$ in 
the second sum, we get
$$
\sum_{0 \leq j < 2^n} \cos((2j+1)x -2^nx) =
2 \sum_{0 \leq t < 2^{n-1}} \cos((2t+1)x)
$$
which gives Identity~\ref{cos2simple}. Similarly Identity~\ref{sin2} implies
Identity~\ref{sin2simple} after noting that for $t \in [0, 2^{2m-1})$ we have
$(-1)^{s(2^{2m-1} - t -1)} = -(-1)^{s(t)}$ and $s(2^{2m-1} + t) = -(-1)^{s(t)}$. \endpf

\begin{remark}\label{epsilon}

\ { }

-- In Identities~\ref{general} to \ref{sin} above, the quantities
$$
\sum_{0 \leq j < 2^n} 
f(s(j)) F\left(\sum_{0 \leq q < n} \sigma_q(j) \lambda_q\right)
$$ 
can be replaced by the (equal) expressions 
$$
\sum_{\varepsilon_0, \varepsilon_1, \cdots, \varepsilon_{n-1}
\in \{-1, +1\}^n} 
f\left(\sum_{0 \leq q < n} \frac{1 + \varepsilon_q}{2}\right)
F\left(\sum_{0 \leq q < n} \varepsilon_q \lambda_q\right).
$$

\bigskip

-- Identities analogous to Identities~\ref{sinh}, \ref{sin}, \ref{sin2}, 
and \ref{sin2simple} above, where $2m$ is replaced by $2m+1$, can
be obtained with similar arguments. For example 
$$
\prod_{0 \leq k < 2m+1} \sin \lambda_k
= (-1)^{m+1}2^{-2m-1} \sum_{0 \leq j < 2^{2m+1}} (-1)^{s(j)}
\sin\left(\sum_{0 \leq q < 2m+1} \sigma_q(j) \lambda_q\right)
$$
(this relation can be deduced directly from Identity~\ref{sin} by 
multiplicating by $\sin \lambda_{2m}$ and using 
$\sin \alpha \cos \beta = (\sin(\alpha + \beta) +\sin(\alpha - \beta))/2$).
\end{remark}

\section*{Wallis formula}

John Wallis (1616-1703), while calculating the value of $\pi$ by finding 
the area under the quadrant of a circle (ten years before Newton discovered 
calculus), obtained the now called ``Wallis formula''. It is worth mentioning
that Sondow showed how to use Wallis formula to compute in an elementary 
way $\zeta(0)$ and $\zeta'(0)$, where $\zeta$ is the Riemann zeta function
\cite{Sondow}. We show here how Wallis formula can be deduced from 
Identity~\ref{sin}.

\begin{corollary}[Wallis formula]
If $n$ tends to infinity, then
$$
\lim\frac{2^{4n}}{n{2n \choose n}^2} = \pi.
$$
This can also be written
$$
\frac{\pi}{2} = \frac{2 \cdot 2}{1 \cdot 3} \ \frac{4 \cdot 4}{3 \cdot 5}
\ \frac{6 \cdot 6}{5 \cdot 7} \ \cdots
$$
\end{corollary}

\proof  Applying Identity~\ref{sin} with 
$\lambda_0 = \lambda_1 = \cdots = \lambda_{2m-1} := x$ we obtain
$$
(\sin x)^{2m} = (-1)^m 2^{-2m} 
\sum_{0 \leq j < 2^{2m}} (-1)^{s(j)} 
\cos\left(x\sum_{0 \leq q < 2m} \sigma_q(j)\right).
$$
Since $\sigma_q(j) = 2e_q(j) - 1$, we have
$\sum_{0 \leq q < 2m} \sigma_q(j) = 2s(j) - 2m$, hence
$$
(\sin x)^{2m} = (-1)^m 2^{-2m}
\sum_{0 \leq j < 2^{2m}} (-1)^{s(j)}
\cos(2x(s(j)-m)).
$$
Integrating yields
$$
\int_0^{\pi/2} (\sin x)^{2m} dx =
(-1)^m 2^{-2m} \sum_{0 \leq j < 2^{2m}} (-1)^{s(j)}
\int_0^{\pi/2} \cos(2x(s(j)-m))dx.
$$
The integrals in the right-hand quantity are equal to 
zero if $s(j) \neq m$ and to $\pi/2$ if $s(j)=m$. 
There are exactly ${2m \choose m}$ integers $j \in [0, 2^{2m})$ for 
which $s(j)=m$, hence
$$
\int_0^{\pi/2} (\sin x)^{2m} dx =
(-1)^m 2^{-2m}(-1)^m {2m \choose m} \frac{\pi}{2} 
= \frac{1}{2^{2m}} {2m \choose m} \frac{\pi}{2}.
$$
Now, using the classical trick, define 
$I_n := \int_0^{\pi/2} (\sin x)^n dx$. Then
$$
\int_0^{\pi/2} (\sin x)^{n+2} dx = 
\int_0^{\pi/2} (\sin x)^n dx - 
\int_0^{\pi/2} \cos x \ (\cos x (\sin x)^n) dx .
$$
Integrating the last integral by parts gives
$$
(n+2) I_{n+2} = (n+1) I_n \ \mbox{\rm hence }
(n+2) I_{n+2} I_{n+1} = (n+1) I_{n+1} I_n.
$$
This shows that the quantity $(n+1) I_{n+1} I_n$ does not
depend from $n$, hence is equal to $I_1 I_0 = \pi/2$.
Furthemore the equality $(n+2) I_{n+2} = (n+1) I_n$ shows
that $I_{n+2}/I_n$ tends to $1$ as $n$ tends to infinity.
But $I_{n+2} \leq I_{n+1} \leq I_n$; hence dividing by $I_n$
shows that $I_{n+1}/I_n$ tends to $1$ as $n$ tends to infinity.
Hence
$$
\frac{\pi}{2} = (n+1) I_{n+1} I_n \sim n I_n^2
$$
which gives
$$
\pi = \lim_{n \to \infty} 2n I_n^2 = \lim_{n \to \infty} 4n I_{2n}^2 
= \lim_{n \to \infty} \frac{n\pi^2}{2^{4n}} {2n \choose n}^2.
$$
Hence
$$
\lim_{n \to \infty} \frac{2^{4n}}{n {2n \choose n}^2} = \pi
$$

\section*{Pisot numbers}

Pisot-Vijayaraghavan numbers were first studied by C.~Pisot (1910-1988) 
and T.~Vijayaraghavan (1902-1955): these are the algebraic integers 
$> 1$ such that all their other algebraic conjugates lie in the open 
disk $z < 1$. Their principal property is to behave almost like 
integers in questions of distribution modulo $1$.

The product formulas above have applications for the distribution
modulo $1$ of quantities involving Pisot-Vijayaraghavan numbers.
Before stating such an application (Corollary~\ref{Pisot1} below)
we need a lemma.

\begin{lemma}\label{prodcos}
Let $(\lambda_q)_{q \geq 0}$ be a sequence of real numbers. For any 
real number $x$ define $\Vert x \Vert := \min_{x \in \mathbb Z} |x-n|$.
Then the infinite product $\prod_{q \geq 0} |\cos \pi \lambda_q|$ 
equals $0$ if and only if either 
$\sum_{q \geq 0} \Vert \lambda_q\Vert^2 = + \infty$, or there exists 
$q \geq 0$ such that $\lambda_q \equiv \frac{1}{2} \bmod 1$.
\end{lemma}

\proof  Left to the reader.  \endpf

\begin{corollary}[Mend\`es France \cite{Mendes}]\label{Pisot1}
Let $\lambda := (\lambda_j)_{j \in {\mathbb N}}$ be a sequence of real
numbers. For every integer $j$ with binary expansion
$j=\sum_q e_q(j) 2^q$, where $e_q(j) \in \{0, 1\}$, define as previously
$u(j, \lambda) := \sum_q e_q(j) \lambda_q$. Then the sequence 
$j \to u(j, \lambda)$ is uniformly distributed modulo $1$ if and only 
if for each $\ell \in {\mathbb N} \setminus \{0\}$ either
$\sum_q \Vert \ell \lambda_q \Vert^2 = + \infty$ or there exists
$q = q(\ell) \geq 0$ such that $\ell \lambda_q \equiv \frac{1}{2} \bmod 1$.

\medskip

As a consequence, supposing $\theta > 1$,  the sequence 
$(x\sum_q e_q(k) \theta^q)_{k \in {\mathbb N}}$ is uniformly 
distributed modulo $1$ if and only if either $\theta$ is not a
Pisot number and $x \neq 0$, or $\theta$ is a Pisot number and 
$x \notin {\mathbb Q}(\theta)$.

\medskip

In particular, if $x \neq 0$ and $\theta$ is not a Pisot number, 
then the set of finite sums $x \sum \theta^{n_j}$ is dense modulo $1$. 
The same conclusion holds if $\theta$ is a Pisot number and 
$x \notin {\mathbb Q}(\theta)$.
\end{corollary}

\proof We first recall that $\theta$ is a Pisot number if and only if
there exists a nonzero real number $x$ such that
$\sum_q \Vert x \theta^q \Vert^2 < \infty$, see \cite[p. 238]{Pisot};
furthermore any such $x$ must necessarily belong to ${\mathbb Q}(\theta)$
see \cite[Th\'eor\`eme 2, p. 153-154]{Pisot2} or \cite[Theorem A, p. 4]{Salem}.
Now, since the sequence $j \to u(j, \lambda)$ is uniformly distributed 
modulo $1$ if and only if 
$$
\lim_{N \to \infty}
\left(\frac{1}{N}\sum_{0 \leq j < N} \exp(2i \pi \ell u(j,\lambda))\right) = 0
$$
for all integers $\ell \neq 0$, we see from Lemma~\ref{prodcos} that it 
suffices to prove that conditions (i) and (ii) below are equivalent:
$$
\begin{array}{lll}
&(i) \ \ &\forall \ell \neq 0, \ \ 
\displaystyle\lim_{N \to \infty}\frac{1}{N}\left(\sum_{0 \leq j < N} 
\exp(2i \pi \ell u(j,\lambda))\right) = 0 \\
&(ii) \ \  
&\forall \ell \neq 0, \ \ \ 
\displaystyle\prod_{q \geq 0} |\cos \pi \ell \lambda_q| = 0.
\end{array}
$$
Now Identity~\ref{step} with $a := 1$ and $Z := 2i\pi\ell$ 
(where $\ell$ is an integer) reads
$$
\prod_{0 \leq k < n} \left(\frac{1 + e^{2i\pi \ell \lambda_k}}{2}\right) 
= \frac{1}{2^n} \sum_{0 \leq j < 2^n} e^{2i\pi \ell u(j, \lambda)}.
$$
Hence
$$
\prod_{0 \leq k < n} |\cos \pi \ell \lambda_k| =
\left|\frac{1}{2^n}\sum_{0 \leq j < 2^n} e^{2i\pi \ell u(j, \lambda)} \right|
$$
showing that (i) implies (ii).

\medskip

Conversely, suppose that Condition (ii) holds. Note that 
$u(2j, \lambda) = u(j, T\lambda)$ and that $u(2j+1, \lambda) = \lambda_0
+ u(j, T\lambda)$, where $T\lambda$ is the shifted sequence
$T\lambda := (\lambda_{q+1})_{q \geq 0}$. Hence, defining ${\cal M}(N, \lambda)$
by 
$$
{\cal M}(\lambda) := 
\limsup_{N \to \infty}
\left|\frac{1}{N} \sum_{0 \leq k < N} \exp(2i \pi \ell u(k, \lambda))\right|
$$
we have (see Lemma~3 of \cite{Mendes})
$$
{\cal M}(\lambda) = 
|\cos(\pi \ell \lambda_0)| {\cal M}(T\lambda) 
= |\cos(\pi \ell \lambda_0)| |\cos(\pi \ell \lambda_1)| {\cal M}(T^2\lambda)  
= \ldots
$$
hence
$$
{\cal M}(\lambda) \leq \prod_{q \geq 0} |\cos(\pi \ell \lambda_q)| = 0.
\hskip 3truecm \Box
$$
\section*{Multigrade equalities and the Prouhet-Tarry-Escott problem}
\label{multigrades}

We alluded at the beginning of this paper to multigrade equalities.
The pioneering work in this area is the very short Note aux 
Comptes-Rendus de l'Acad\'emie des Sciences of Prouhet in 1851: for 
a survey on the Prouhet problem see for example \cite{BorIng, Salon}. 
In this section we sketch the construction of multigrade equalities \`a 
la Prouhet using our fundamental identity.

\bigskip

Replacing $\lambda_k$ by $x \lambda_k$ in Identity~\ref{sinh} gives
$$
\begin{array}{lll}
\displaystyle\prod_{0 \leq k < 2m} \sinh x \lambda_k
&=& 2^{-2m} \displaystyle\sum_{0 \leq j < 2^{2m}} (-1)^{s(j)}
\cosh\left(x \sum_{0 \leq q < 2m} \sigma_q(j) \lambda_q\right) \\
&=& 2^{-2m} \displaystyle\sum_{0 \leq j < 2^{2m}} (-1)^{s(j)}
\sum_{\ell \geq 0} \frac{1}{(2\ell)!} x^{2\ell}
\left(\sum_{0 \leq q < 2m} \sigma_q(j) \lambda_q\right)^{2\ell} \\
&=& 2^{-2m} \displaystyle\sum_{\ell \geq 0} \frac{1}{(2\ell)!} x^{2\ell}
\sum_{0 \leq j < 2^{2m}} (-1)^{s(j)}
\left(\sum_{0 \leq q < 2m} \sigma_q(j) \lambda_q\right)^{2\ell}.
\end{array}
$$
Since the left-hand expression is equal to 
$\lambda_0 \lambda_1 \cdots \lambda_{2m-1} x^{2m} + {\cal O}(x^{2m+2})$ 
in the neighborhood of $0$, we have
for $\ell \in [0, m-1]$
$$
\sum_{0 \leq j < 2^{2m}} (-1)^{s(j)}
\left[\left(\sum_{0 \leq q < 2m} \sigma_q(j) \lambda_q\right)^2\right]^\ell
= 0.
$$
Choosing the $\lambda_j$'s adequately gives a family of multigrade
equalities.

\bigskip

In the same way, using the second part of Remark~\ref{epsilon}, one
can obtain an identity for $\prod_{0 \leq k < 2m+1} \sinh x \lambda_k$,
and hence another family of multigrade identities.

\section*{Back to iterating the duplication of sines}

\subsection*{A natural occurrence of morphisms of the free monoid}

We recall the notation: 
$P_n(x) := \sin x \sin 2x \sin 4x \cdots \sin 2^n x$.
We will denote by ${\rm sgn\,} y$ the sign of the real number $y$,
i.e., ${\rm sgn\,} y := y/|y|$ if $y \neq 0$. We will show that the
signs of $P_n(x)$, for $x \in (0, \pi)$, form a word of length $2^n$
that is obtained by iterating the morphism $\varphi$
of the free monoid generated by $\{0, 1\}$ defined by
$$
\begin{array}{lll}
\varphi(+) &=& + - \\
\varphi(-) &=& - + \\
\end{array}
$$
(in particular when $n$ goes to infinity, the sequence of signs of
$P_n(x)$ converges to the Prouhet-Thue-Morse sequence).

\begin{proposition}
The finite sequence of (consecutive distinct) signs of $P_{n+1}(x)$ for 
$x \in (0, \pi)$ can be obtained from the sequence of signs of $P_n(x)$ 
for $x \in (0, \pi)$ by applying the morphism of monoid
$+ \to + \ -$, $- \to - \ +$. In particular the sequence of signs
of $P_n(x)$ for $x \in (0, \pi)$ is the prefix of length $2^n$
of the Prouhet-Thue-Morse sequence.
\end{proposition}

\proof For $n \geq 0$ and $j \in [0, 2^n)$, let $p_{n,j}$ be the sign
of $P_n(x)$ for $x \in (\frac{j\pi}{2^n}, \frac{(j+1)\pi}{2^n})$
(it is clear that $P_n(x)$ has constant sign on such an interval).

Let now $j$ belong to $[0, 2^{n+1})$ and let us determine the sign of 
$P_{n+1}(x)$ for $x \in (\frac{j\pi}{2^{n+1}}, \frac{(j+1)\pi}{2^{n+1}})$.

\begin{itemize}

\item If $j$ is even, say $j := 2k$, the interval 
$(\frac{j\pi}{2^{n+1}}, \frac{(j+1)\pi}{2^{n+1}})$ is
contained in the interval 
$(\frac{u\pi}{2^n}, \frac{(u+1)\pi}{2^n})$
and $2^{n+1} x$ belongs to $(2u\pi, 2(u+1)\pi$, hence
$p_{n+1,j} = p_{n+1,2u} = p_{n,u}$.

\item If $j$ is odd, say $j := 2k$, the interval 
$(\frac{j\pi}{2^{n+1}}, \frac{(j+1)\pi}{2^{n+1}})$ is
contained in the interval 
$(\frac{u\pi}{2^n}, \frac{(u+1)\pi}{2^n})$
and $2^{n+1} x$ belongs to $((2u+1)\pi, (2u+2)\pi$, hence
$p_{n+1,j} = p_{n+1,2u+1} = -p_{n,u}$.

\end{itemize}

The formulas giving $p_{n+1,j}$ in terms of $p_{n,j}$ exactly mean
that each of the intervals $(\frac{j\pi}{2^n}, \frac{(j+1)\pi}{2^n})$
where $P_n(x)$ has constant sign splits into the two subintervals 
$(\frac{2j\pi}{2^{n+1}}, \frac{(2j+1)\pi}{2^{n+1}})$ and
$(\frac{(2j+1)\pi}{2^{n+1}}, \frac{(2j+2)\pi}{2^{n+1}})$
where $P_{n+1}(x)$ has constant sign, this sign being the same
as the sign of $P_n(x)$ on the first interval and the opposite on
the second one. This can be described by 
$$
\begin{array}{lll}
+ &\to& + \ - \\
- &\to& - \ + 
\end{array}
$$
at each splitting.
\endpf

\begin{remark}
What seems interesting in this (easy) result is the natural
occurrence of a morphism of monoid in a simple question that
is not explicitly related to combinatorics on words.
For more material on morphisms of (free) monoids, the reader might
want to read for example \cite{AllSha}.
\end{remark}

\subsection*{Asymptotic behavior of the duplicating sinus product}

In this section we will study the asymptotic behavior of the product
$$
P_n(x) := \sin x \sin 2x \sin 4x \cdots \sin 2^n x,
$$
in particular of the quantities
$\Vert P_n \Vert_{\infty} := \sup_{x \in {\mathbb R}} |P_n(x)|$
and
$\Vert P_n \Vert_1 := \int_0^{\pi} |P_n(x)| dx$.

\bigskip

\begin{theorem}[Gelfond \cite{Gelfond}, Newman-Slater \cite{Newman-Slater}]
\label{Gelfond}
The following bound holds
$$
\Vert P_n \Vert_{\infty} \leq \left(\frac{\sqrt 3}{2}\right)^n.
$$
\end{theorem}

\begin{remark}

\ { }

-- Note that this inequality is nearly optimal 
since $|P_n(\pi/3)| = ({\sqrt 3}/2)^{n+1}$. 

\medskip

-- This result seems to have been proved for the first time by
Gelfond who used a simple and ingenious lemma \cite[p. 62]{Gelfond}.
Gelfond applied this bound to obtain an asymptotic estimation of the
number of integers $\leq x$ belonging to a given class of congruence
and whose sum of digits in some integral base belongs to a given
class of congruence. Another proof was given by Newman and Slater in
\cite[p. 73]{Newman-Slater}: they used this result to establish the 
preponderance of integers having an even number of $1$'s in their binary 
expansion over integers having an odd number of $1$'s in natural sequences 
of integers (for example the sequence $(M_n)_{n \geq 0}$, or the sequence 
of squarefree integers). Among other occurrences of this result in the 
literature we cite the paper of Fouvry and Mauduit 
\cite[p. 583]{Fouvry-Mauduit} where they were interested {\it inter alia\,}
in almost-primes in sequences of integers with automatic characteristic 
function; their proof is the same as the proof of Newman and Slater in
\cite{Newman-Slater}.

Note that $\log(({\sqrt 3}/2)^n) = n({\log 3})/({2\log 2})$ and that 
the exponent $\frac{\log 3}{2 \log 2}$ as well as the summatory function of the 
sequence $((-1)^{s(3n)})_{n \geq 0}$ also occur in \cite{Newman, Coquet}.

\end{remark}

For the sake of completeness we give below the proof of Newman and Slater
slightly generalized to the case where duplication is replaced by
multiplication by any integer. (For recovering Theorem~\ref{Gelfond} above
take $r := 2$ in Proposition~\ref{r-Gelfond} below.)

\begin{proposition}\label{r-Gelfond}
Let $r$ be a positive integer. Define $P_{r,n}(x)$ by
$$
P_{r,n}(x) := |\sin x \sin rx \cdots \sin r^n x|.
$$

If $r$ is odd, then \ \ \ 
$\Vert P_{r,n} \Vert_{\infty} = 1$,

\bigskip

If $r$ is even, then \ \ \  
$\displaystyle
\left(\cos\left(\frac{\pi}{2r+2}\right)\right)^{n+1} 
\leq \Vert P_{r,n} \Vert_{\infty} \leq
\left(\cos\left(\frac{\pi}{2r+2}\right)\right)^n$.
\end{proposition}

\proof To avoid both absolute values and taking real powers of negative 
numbers, we will study the quantity $P_{r,n}^2(x)$. The following trick
is adapted from \cite{Newman-Slater}:
$$
P_{r,n}^2(x) = (\sin^2 x)^{1/r+1} 
\left(\prod_{0 \leq j \leq n-1} (\sin^2 r^j x)^{r/r+1} (\sin^2 r^{j+1} x)^{1/r+1}\right)
(\sin^2 r^n x)^{r/r+1}.
$$
Hence, defining the function $\varphi$ by 
$$
\varphi(x) := (\sin^2 x)^{r/r+1} (\sin^2 r x)^{1/r+1},
$$
we have
$$
P_{r,n}^2(x) \leq \prod_{0 \leq j \leq n-1} \varphi(r^j x).
$$
Now the derivative of the function $\varphi$ is given by:
$$
\varphi' (x) = \frac{2r}{r+1} (\sin^2 x)^{-1/r+1} (\sin^2 r x)^{-r/r+1}
\sin x \sin rx \sin((r+1)x).
$$
In particular the maximum of $\varphi$ must occur for $x$ such that
$\sin((r+1)x) = 0$ i.e., $x \in \{ k \pi /(r+1), \ k \geq 0 \}$.
It is thus immediate that the maximum of $\varphi$ occurs at
$x = \frac{\lfloor\frac{r+1}{2}\rfloor}{r+1}\pi$ and that this maximum
is equal to $\sin^2 \left(\frac{\lfloor\frac{r+1}{2}\rfloor}{r+1}\pi\right)$.
Hence 
$$
P_{r,n}^2(x) \leq 
\left(\sin^2 \left(\frac{\lfloor\frac{r+1}{2}\rfloor}{r+1}\pi\right)\right)^n
$$
which yields an upper bound for $\Vert P_{r,n} \Vert_{\infty}$.

\bigskip

Taking $x_0 := (\pi\lfloor\frac{r+1}{2}\rfloor)/(r+1)$, we see that, 
for any integer $k \geq 0$,
$$
r^{k+1}x_0 + r^k x_0 = r^k \pi \left\lfloor\frac{r+1}{2}\right\rfloor 
\in \pi {\mathbb Z}^+
$$
hence
$$
|\sin r^{k+1} x_0| = |\sin r^k x_0| = \ldots = |\sin x_0|.
$$
This implies that $P_{r,n}(x_0) = |\sin x_0|^{n+1}$, hence
$$
\left(\sin 
\left(\frac{\lfloor\frac{r+1}{2}\rfloor}{r+1}\pi\right)\right)^{n+1}
\leq \Vert P_{r,n} \Vert_{\infty} \leq
\left(\sin
\left(\frac{\lfloor\frac{r+1}{2}\rfloor}{r+1}\pi\right)\right)^n
$$
which gives the desired result.  \endpf

\begin{remark}
Proposition~\ref{r-Gelfond} shows in particular that, for all $n$,
$\lim_{r \to \infty} \Vert P_{n,r} \Vert_{\infty} = 1$.
\end{remark}

\bigskip

The asymptotic behavior of the $L_1$-norm of $P_n(x) := P_{2,n}(x)$ also 
occurs in the literature. We have the following result.

\begin{theorem}[Fouvry-Mauduit \cite{Fouvry-Mauduit}]
\label{Fo-Ma}
The following equality holds
$$
\Vert P_n \Vert_1 = \int_0^{\pi} |P_n(t)| dt = C \rho^n (1 + o(1))
$$
where $C$ is a positive constant and $\rho$ satisfies
$0.654336 \cdots < \rho < 0.663197 \cdots$.
\end{theorem}

\begin{remark}

Theorem~\ref{Fo-Ma} above is proved in \cite{Fouvry-Mauduit} where
the authors are interested in almost-primes (i.e., integers having
at most a given number of primes) in sequences of integers whose
characteristic functions are automatic; their proof uses transfer
operators. Note that \`Eminyan gave earlier in \cite{Eminyan1} a bound
for $\Vert P_{2n-1} \Vert_1$ of the form $C' \mu^{2n-1}$ where
$\mu = (2 + {\sqrt 2})^{1/4} /2$. This result is used to compute 
the average value of the number of divisors of $n$ for $n$ having
an even binary sum of digits. This result is also used in 
\cite{Eminyan2} to compute the average value of the number of
representations of $n$ as a sum of two squares for $n$ having an
even binary sum of digits. 
\end{remark}

A rough upper bound for $\rho$ in Theorem~\ref{Fo-Ma} can easily be 
obtained from Identity~\ref{sin2simple} and Parseval's equality 
(orthogonality of sines or cosines). We state this result in a more 
general framework.

\begin{proposition}\label{lambda-epsilon}
Let $\lambda = (\lambda_0, \ldots, \lambda_N)$ be an $(N+1)$-tuple
of integers such that, for all choices of $\varepsilon_q = \pm 1$,
and  $\varepsilon_q' = \pm 1$, we have
$$
\sum_{0 \leq q \leq N} \varepsilon_q \lambda_q  =
\sum_{0 \leq q \leq N} \varepsilon_q' \lambda_q
\ \ \ \Rightarrow \  \ \varepsilon_q = \varepsilon_q' \ 
\mbox{\rm for all} \ q.
$$
Define
$$
P_N^{\lambda}(x) := \left|\prod_{0 \leq j \leq N} \sin(\lambda_j x)\right|.
$$
Then
$$
\int_0^{\pi} |P_N^{\lambda}(x)| dx \leq \frac{\pi}{2^{1+N/2}}.
$$
\end{proposition}

\begin{remark}
If $\lambda_q := 2^q$ for all $q$, the condition on $\varepsilon_q$
and $\varepsilon_q'$ holds. Proposition~\ref{lambda-epsilon} yields for $\rho$ 
in Theorem~\ref{Fo-Ma} the bound $\rho \leq \frac{\sqrt{2}}{2} = 0.707106 \cdots$.
This bound is of course less precise than the result of \cite{Fouvry-Mauduit}, but 
more precise than the bound obtained from the inequality 
$\Vert P_n \Vert_1 \leq C \Vert P_n \Vert_{\infty}$ which 
gives $\rho \leq \frac{\sqrt{3}}{2} = 0.866025 \cdots$.
\end{remark}

\proof The Fourier expansion of $P_N^{\lambda}$ depends on the parity
of $N$: we have from Identity~\ref{sin} and the two parts of
Remark~\ref{epsilon}, where $\lambda_k$ is replaced by $\lambda_k x$:
$$
P_N^{\lambda}(x) = \frac{1}{2^{N+1}} \left|
\sum_{\varepsilon \in \{-1, +1\}^{N+1}} 
\pm F\left(x \sum_{0 \leq q \leq N} \varepsilon_q \lambda_q\right) \right|
$$
where the $\pm$ symbol equals $\varepsilon_0 \varepsilon_1 \cdots \varepsilon_N$,
and where $F(\xi) = \cos \xi$ if $N$ is odd and $F(\xi) = \sin \xi$ if $N$ 
is even. In either case, Parseval's equality reads
$$
\frac{1}{\pi} \int_0^{\pi} |P_N^{\lambda}(x)|^2 dx
=\frac{1}{2} \sum_{0 \leq k \leq 2^{N+1} - 1} \frac{1}{2^{2(N+1)}}
= \frac{1}{2^{N+2}}.
$$
Applying the inequality of Cauchy-Schwarz yields
$$
\frac{1}{\pi} \int_0^{\pi} |P_N^{\lambda}(x)|dx \leq
\left(\frac{1}{\pi} \int_0^{\pi} |P_N^{\lambda}(x)|^2 dx\right)^{1/2}
= \frac{1}{2^{1+N/2}}.
\hskip 3truecm \Box
$$
A last comment. The sup norm of the product $|\prod\cos(\lambda_j x)|$ 
is obviously $1$. Its $L^2$-norm can be computed as above:
the graph of the product of sines or cosines $L^2$-collapses onto the
horizontal axis as $N$ increases.

\end{document}